# Riemann problem for a nonsymmetric Keyfitz-Kranzer and pressureless gas systems with a time-dependent Coulomb-like friction term

Richard De la cruz[a,*], Wladimir Neves[b]

[a]*School of Mathematics and Statistics, Universidad Pedagógica y Tecnológica de Colombia, 150003, Colombia.*
[b]*Instituto de Matemática, Universidade Federal do Rio de Janeiro, Cidade Universitária 21945-970, Rio de Janeiro, Brazil.*

**ARTICLE INFO**

**ABSTRACT**

In this paper, we study the Riemann solutions for two systems: the nonsymmetric Keyfitz-Kranzer system and the pressureless system, both of which have a time-dependent Coulomb-like friction term. Our analysis identified two types of Riemann solutions: contact discontinuities and delta-shock solutions. We obtain generalized Rankine-Hugoniot conditions, which is the support for constructing the delta-shock solution for the nonsymmetric Keyfitz-Kranzer system with a time-dependent Coulomb-like friction term. Furthermore, we demonstrate that as the pressure tends to zero, the Riemann solutions of the nonsymmetric Keyfitz-Kranzer system converge to those of the pressureless system, with both systems incorporating a time-dependent Coulomb-like friction term.

## 1. Introduction

In this paper, we study the Riemann solutions for the following two systems:

$$\begin{cases} \rho_t + (\rho u)_x = 0, \\ (\rho u)_t + (\rho u(u + \frac{\mu}{\rho}))_x = \alpha(t)\rho \end{cases} \tag{1}$$

and

$$\begin{cases} \rho_t + (\rho u)_x = 0, \\ (\rho u)_t + (\rho u^2)_x = \alpha(t)\rho, \end{cases} \tag{2}$$

where $\mu > 0$ and $\alpha \in C(\mathbb{R})$. The system (1) is a prototype of the nonsymmetric Keyfitz-Kranzer system whereas the system (2) is the pressureless gas system. In both systems, the source term corresponds to the Coulomb-like friction term, which changes with time. To study the Riemann problem for the systems (1) and (2), we consider the following initial data:

$$(\rho(x,0), u(x,0)) = \begin{cases} (\rho_-, u_-), & \text{if } x < 0, \\ (\rho_+, u_+), & \text{if } x > 0, \end{cases} \tag{3}$$

for arbitrary constant states $(\rho_\pm, u_\pm)$ with $\rho_\pm > 0$. These systems are prevalent in physics and engineering models, prompting extensive study by numerous authors and resulting in several rigorous findings. In both systems (1) and (2), the homogeneous case arises when $\alpha(\cdot) \equiv 0$.

We note that the homogeneous case of system (1) is a particular case of the following system

$$\begin{cases} \rho_t + (\rho\Psi(\rho, u))_x = 0, \\ (\rho u)_t + (\rho u\Psi(\rho, u))_x = 0, \end{cases} \tag{4}$$

*Corresponding author
✉ richard.delacruz@uptc.edu.co (R. De la cruz); wladimir@im.ufrj.br (W. Neves)
ORCID(s): 0000-0003-1342-7946 (R. De la cruz); 0000-0000-0000-0000 (W. Neves)

which finds applications in various physical contexts contingent on the function $\Psi$. For example, if we replace $\rho$ by $v$ and $\rho u$ by $w$, the system (4) can be rewritten into a system in the variables $v$ and $w$ which was proposed, by Keyfitz and Kranzer [19], as a model for the elastic string. From the remarkable work of Keyfitz and Kranzer [19], this system was called the *Keyfitz-Kranzer system*.

Additionally, system (4) serves as a model for polymer flooding in oil reservoirs [17]. Temple [30] investigated the system (4) under the condition that $\rho\Psi(\rho, u)$ is not a convex function. In this scenario, the system describes how the introduction of a polymer affects the flow of water and oil within a reservoir. Temple established the existence of a global weak solution to the Cauchy problem. In contrast, Aw and Rascle [2] examined the system (4) with $\Psi(\rho, u) = u - P(\rho)$, proposing it as a macroscopic model for traffic flow, where $\rho$ and $u$ represent the density and velocity of cars on a roadway, respectively, and $P$ is a smooth and strictly increasing function satisfying $\rho P''(\rho) + 2P'(\rho) > 0$ for $\rho > 0$. In a broader context, Lu [23] demonstrated the existence of a global weak solution to the Cauchy problem for the system (4) with $\Psi(\rho, u) = f(u) - P(\rho)$ where $f$ is a non-negative convex function and $P$ satisfies the following conditions:

$$P(\rho) \le 0 \text{ for } \rho > 0, \ P(0) = 0, \ \lim_{\rho \to 0} \rho P'(\rho) = 0, \ \lim_{\rho \to \infty} P(\rho) = \infty, \text{ and } \rho P''(\rho) + 2P'(\rho) < 0 \text{ for } \rho > 0.$$

De la cruz and Santos [10] addressed the Riemann problem of the system (4) with $\Psi(\rho, u) = f(u) - P(\rho)$, where $f$ is a strictly increasing function and $P(\rho) = -\frac{1}{\rho}$ for $\rho > 0$. The pressure function given by $P(\rho) = -\frac{\mu}{\rho}$ (where $\mu > 0$) is known as *Chaplygin pressure*. Chaplygin introduced this model for pressure as a suitable mathematical approximation for calculating the lifting force on an airplane wing in aerodynamics [7]. We observe that the Chaplygin pressure $P(\rho) = -\frac{\mu}{\rho}$ (where $\mu > 0$) is the prototype function satisfying $\rho P''(\rho) + 2P'(\rho) = 0$ for $\rho > 0$ and it is not possible to solve the Riemann problem for (4) using only classical elementary waves. Consequently, in [8, 10], the existence and uniqueness of Riemann solutions, including delta-shock solutions, were demonstrated. Furthermore, a comprehensive examination of the interaction between delta-shock waves and elementary waves is presented in [11].

On the other hand, the homogeneous case of system (2) is known as the *pressureless gas system* or *one-dimensional zero-pressure gas dynamics system*, and it has garnered extensive attention because of its significance in various applications (see [3, 5, 6, 13, 16, 21] and references therein for further details). The pressureless gas system, also known as the adhesion particle dynamics system, describes the motion process of free particles sticking together upon collision at low temperatures and the formation of large-scale structures in the universe [6, 27]. In particular, the pressureless gas system emerges in a diverse array of physics models (refer to [3, 13, 21, 24] for more information). It is well established that both delta-shock waves and vacuum states occur in the Riemann problem of the pressureless gas system (see [25, 29, 31]). A delta-shock wave solution is a non-classical wave solution in which at least one state variable may develop a weighted Dirac delta measure. In the context of pressureless gas dynamics system, delta-shocks represent phenomena such as galaxy formation or particle concentrations [22].

It is important to note that, the homogeneous systems corresponding to systems (1) and (2) remain invariant under uniform stretching of coordinates: $(x, t) \mapsto (\beta x, \beta t)$ with $\beta$ a constant. Therefore, they admit self-similar solutions defined on the space-time plane, which remains constant along straight-line rays emanating from the origin. However, the structure of solutions for the Riemann problem in inhomogeneous systems (1) and (2) is more complex because of the absence of self-similar solutions in the form of $(\rho(x, t), u(x, t)) = (\rho(x/t), u(x/t))$ due to their inhomogeneity.

One observes that the Coulomb-like friction term $\alpha\rho$ ($\alpha$ being constant) was introduced by Savage and Hutter [26] to describe the granular flow behavior. Abreu, De la cruz, and Lambert [1] showed the existence of Riemann solutions for a nonsymmetric Keyfitz-Kranzer system with Coulomb-like friction, i.e., Riemann solutions for the system (1) with $\alpha(\cdot) \equiv \alpha = const.$ and $\mu = 1$. Specifically, they established the existence of both shock wave solutions (satisfying the classical Rankine-Hugoniot conditions) and intricate delta-shock wave solutions (satisfying the generalized Rankine-Hugoniot conditions). In contrast, Shen [28] delved into the Riemann problem for the pressureless gas system with a Coulomb-like friction term (with $\alpha(\cdot) \equiv \alpha = const.$, i.e., the coefficient remains constant on time), finding solutions consisting of delta-shock waves and vacuum states. Shen's work represents a pioneering effort in tackling the Riemann problem for pressureless gas dynamics systems with a source term. Subsequently, Keita and Bourgault [18] solved the Riemann problem for the pressureless system with linear damping, presenting results that encompass delta-shock wave solutions. Yang [31] solved the Riemann problem for a generalized pressure system, where the Riemann solutions precisely comprise two types: delta-shock waves and vacuum states. Sarrico [25] explored delta-shock solutions for a generalized pressureless gas dynamics system, employing the concept of $\alpha$-solution within the product of distributions. Notably, both studies considered systems without a source term (i.e., homogeneous systems). De la cruz and Juajibioy [12] obtained delta-shock solutions for a generalized pressureless system with linear damping.

Concerning systems (1) and (2), it is natural to inquire about the behavior when $\mu$ tends to 0 in system (1). Specifically, we have the following question

**(Q) :** if $\mu$ approaches 0, do the Riemann solutions of the nonsymmetric Keyfitz-Kranzer system (1) converge to those of the pressureless gas system (2)?

This paper investigates the Riemann problems associated with the nonsymmetric Keyfitz-Kranzer system (1) and the pressureless system (2), both of which feature a time-dependent Coulomb-like friction term. We derive Riemann solutions characterized by contact discontinuities and delta-shock solutions for both systems by employing state variable and similarity techniques. We establish generalized Rankine-Hugoniot conditions to obtain a delta-shock solution under an entropy condition. We analyze the behavior when $\mu$ tends to 0 in system (1) and conclude that Riemann solutions of the nonsymmetric Keyfitz-Kranzer system (1) converge to those of the pressureless gas system (2). Furthermore, we provide an affirmative answer to question (Q).

## 2. The Keyfitz-Kranzer type system

In this section, we study the Riemann problem for (1), that is, the nonsymmetric Keyfitz-Kranzer system with time-dependent Coulomb-like friction term, with initial data (3). For $u_- < u_+ + \mu/\rho_+$, we apply a suitable transformation to obtain an auxiliary system related to the system (1), and thus classical Riemann solutions are obtained. Then, for $u_- > u_+ + \mu/\rho_+$, we constructed a delta-shock solution to the Riemann problem for the nonsymmetric Keyfitz-Kranzer system with a time-dependent Coulomb-like friction term. Here, we obtain generalized Rankine-Hugoniot conditions, which are the background for obtaining the delta-shock solution to the system (1).

### 2.1. Classical Riemann solutions

Motivated by Faccononi and Mangeney [14] for the shallow water equations, we introduce a state variable $\widetilde{u}(x, t)$ to perform the transformation

$$\widetilde{u}(x, t) = u(x, t) - \int_0^t \alpha(s) \, ds. \tag{5}$$

This type of transformation for variables is an effective approach to studying balance laws with external force terms like in (1). Indeed, under the transformation (5) the system (1) and the initial data (3) reduce respectively to

$$\begin{cases} \rho_t + \left( \rho(\widetilde{u} + \int_0^t \alpha(s)ds) \right)_x = 0, \\ (\rho\widetilde{u})_t + \left( \rho\widetilde{u}(\widetilde{u} + \int_0^t \alpha(s)ds) + \mu\widetilde{u} \right)_x = 0, \end{cases} \tag{6}$$

and

$$(\rho(x, 0), \widetilde{u}(x, 0)) = \begin{cases} (\rho_-, u_-), & \text{if } x < 0, \\ (\rho_+, u_+), & \text{if } x > 0. \end{cases} \tag{7}$$

The eigenvalues of the system (6), with respective right eigenvectors, are:

$$\lambda_1(\rho, \widetilde{u}) = \widetilde{u} + \int_0^t \alpha(s) \, ds, \quad \mathbf{r}_1 = (1, 0) \quad \text{and} \quad \lambda_2(\rho, \widetilde{u}) = \widetilde{u} + \frac{\mu}{\rho} + \int_0^t \alpha(s) \, ds, \quad \mathbf{r}_2 = (1, \mu/\rho^2).$$

Observe that both characteristics are linearly degenerate; thus, the system is of the Temple class [30]. Thus, shock wave curves in the phase plane $(\rho, u)$ correspond to the so-called contact discontinuities in the $(x, t)$ plane, which are regular outside characteristic curves, with a jump through it. Moreover, since (6) is a $2 \times 2$ system we have two families of such curves in the $(\rho, u)$ plane. These families through a point $(\rho_-, u_-)$ are:

Family 1: horizontal line $\widetilde{u} = u_-$, $\rho > 0$; Family 2: $(\rho, \widetilde{u})$ such that $\widetilde{u} + \frac{\mu}{\rho} + \int_0^t \alpha(s)ds = u_- + \frac{\mu}{\rho_-} + \int_0^t \alpha(s)ds$, $\rho > 0$.

We denote these families by $J_1 \equiv J_1(\rho_-, u_-)$ and $J_2 \equiv J_2(\rho_-, u_-)$. If $x'(t)$ is the speed of the bounded discontinuity $x = x(t)$, the Rankine-Hugoniot conditions for (6) are

$$x'(t)[\rho] + [\rho(\widetilde{u} + \int_0^t \alpha(s)ds)] = 0 \quad \text{and} \quad x'(t)[\rho\widetilde{u}] + [\rho\widetilde{u}(\widetilde{u} + \int_0^t \alpha(s)ds) + \mu u] = 0 \tag{8}$$

where $[\rho] = \rho_- - \rho_+$, $[\rho\widetilde{u}] = \rho_- u_- - \rho_+ u_+$, $[\rho(\widetilde{u} + \int_0^t \alpha(s)ds)] = \rho_-(u_- + \int_0^t \alpha(s)ds) - \rho_+(u_+ + \int_0^t \alpha(s)ds)$, and $[\rho\widetilde{u}(\widetilde{u} + \int_0^t \alpha(s)ds) + \widetilde{u}] = \rho_- u_-(u_- + \int_0^t \alpha(s)ds) + \mu u_- - \rho_+ u_+(u_+ + \int_0^t \alpha(s)ds) - \mu u_+$. In general, $[q] = q_- - q_+$ is the jump across the discontinuity with $q$. From the Rankine-Hugoniot conditions (8), we find that two constant states $(\rho_-, u_-)$ and $(\rho_+, u_+)$ can be connected by contact discontinuities if and only if

$$J_1 \equiv \begin{cases} x'(t) = \widetilde{u} + \int_0^t \alpha(s)ds, \\ u_- = u_+, \end{cases} \qquad \text{and} \qquad J_2 \equiv \begin{cases} x'(t) = \widetilde{u} + \frac{\mu}{\rho} + \int_0^t \alpha(s)ds, \\ u_- + \frac{\mu}{\rho_-} + \int_0^t \alpha(s)ds = u_+ + \frac{\mu}{\rho_+} + \int_0^t \alpha(s)ds. \end{cases}$$

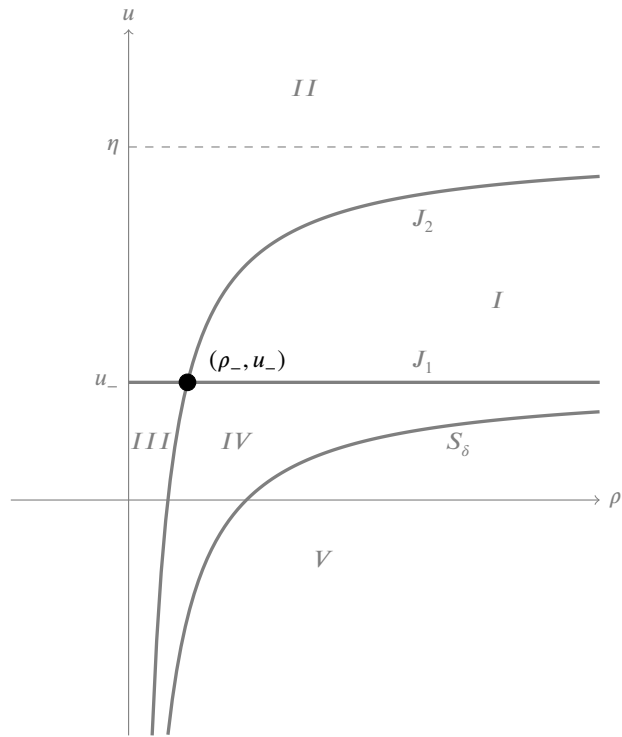

**Figure 1:** Curves $J_1$ and $J_2$ on the phase plane $\rho - \widetilde{u}$. The phase plane is divided into five regions for the constant state $(\rho_-, u_-)$.

For the left state $(\rho_-, u_-)$, the phase plane is divided into five regions (see Figure 1):

1. $I = \{(\rho, \widetilde{u}) : u_- < \widetilde{u} < \eta, \ \mu/v_1 < \rho < +\infty\}$,
2. $II = \{(\rho, \widetilde{u}) : u_- < \widetilde{u} < \eta, \ 0 < \rho < \mu/v_1 \text{ and } \eta \leq \widetilde{u} < +\infty, \ 0 < \rho < +\infty\}$,
3. $III = \{(\rho, \widetilde{u}) : -\infty < \widetilde{u} < u_-, \ 0 < \rho < \mu/v_1\}$,
4. $IV = \{(\rho, \widetilde{u}) : -\infty < \widetilde{u} < u_-, \ \mu/v_1 < \rho < \mu/v_2\}$,
5. $V = \{(\rho, \widetilde{u}) : -\infty < \widetilde{u} < u_-, \ \mu/v_2 < \rho < +\infty\}$,

where $\eta = u_- + \frac{\mu}{\rho_-}$, $v_1 = u_- + \frac{\mu}{\rho_-} - \widetilde{u}$, and $v_2 = u_- - \widetilde{u}$. In the next theorem, we obtain the Riemann solution for $(\rho_+, u_+) \in I \cup II \cup III \cup IV$.

**Theorem 1.** *Let the left and right constant states $(\rho_-, u_-)$ and $(\rho_+, u_+)$ be in $I \cup II \cup III \cup IV$, respectively, such that $u_- < u_+ + \mu/\rho_+$. Then, the Riemann problem* (6) *admits one and only one global solution in the form*

$$(\rho(x,t), \widetilde{u}(x,t)) = \begin{cases} (\rho_-, u_-), & \text{if } x < x_1(t), \\ (\rho_*, u_*), & \text{if } x_1(t) < x < x_2(t), \\ (\rho_+, u_+), & \text{if } x > x_2(t), \end{cases} \tag{9}$$

*where*

$$u_* = u_- \quad and \quad u_* + \frac{\mu}{\rho_*} + \int_0^t \alpha(s)\,ds = u_+ + \frac{\mu}{\rho_+} + \int_0^t \alpha(s)\,ds, \tag{10}$$

*and the contact discontinuities are, respectively, given by*

$$x_1(t) = u_- t + \int_0^t \int_0^r \alpha(s)\,ds\,dr \quad and \quad x_2(t) = \left(u_+ + \frac{\mu}{\rho_+}\right)t + \int_0^t \int_0^r \alpha(s)\,ds\,dr. \tag{11}$$

PROOF. For fixed $(\rho_-, u_-)$, we consider the family of curves

$$F = \{J_2(\rho_+, u_*) \ : \ (\rho_+, u_*) \in J_1(\rho_-, u_-)\}.$$

For $(\rho_+, u_+) \in I \cup II \cup III \cup IV$, the Riemann solution consists of two different (or just one) contact discontinuities. First, there is a 1-contact discontinuity, denoted by $J_1$, connecting $(\rho_-, u_-)$ to intermediate state $(\rho_*, u_*)$, which satisfies Eq. (10). From $(\rho_*, u_*)$, we obtain a 2-contact discontinuity, denoted as $J_2$, to the state $(\rho_+, u_+)$. As we have

$$u_- + \int_0^t \alpha(s)\,ds < u_+ + \frac{\mu}{\rho_+} + \int_0^t \alpha(s)\,ds$$

for all $t \geq 0$, then the solution for the Riemann problem (6)-(7) is given by (9)-(11) and we have $x_1(t) < x_2(t)$ for all $t > 0$.

The solution of the problem (1)-(3) when $u_- < u_+ + \mu/\rho_+$, is directly obtained from the corresponding ones to (6) by performing the transformation of state variables $(\rho(x,t), u(x,t)) = (\rho(x,t), \widetilde{u}(x,t) + \int_0^t \alpha(s)\,ds)$ (see transformation (5)), in which the positions of the contact discontinuities remain unchanged. This fact leads to the following

**Theorem 2.** *Let the left and right constant states $(\rho_-, u_-)$ and $(\rho_+, u_+)$ be in $I \cup II \cup III \cup IV$, respectively, such that $u_- < u_+ + \mu/\rho_+$. Then, the Riemann problem* (1)-(3) *with time-dependent Coulomb-like friction term admits one and only one global solution satisfying* (10) *and* (11) *of the form*

$$(\rho(x,t), u(x,t)) = \begin{cases} (\rho_-, u_- + \int_0^t \alpha(s)\,ds), & \text{if } x < x_1(t), \\ (\rho_*, u_* + \int_0^t \alpha(s)\,ds), & \text{if } x_1(t) < x < x_2(t), \\ (\rho_+, u_+ + \int_0^t \alpha(s)\,ds), & \text{if } x > x_2(t), \end{cases} \tag{12}$$

*with $u_* = u_-$ and $\rho_* = \frac{\mu \rho_+}{\mu + \rho_+(u_+ - u_-)}$.*

We observe that the solution obtained in the previous theorem satisfies the following Lax shock condition

$$\lambda_1(\rho_-, u_- + \int_0^t \alpha(s)\,ds) = \lambda_1(\rho_*, u_* + \int_0^t \alpha(s)\,ds), \quad \frac{dx_1(t)}{dt} < \lambda_2(\rho_*, u_* + \int_0^t \alpha(s)\,ds)$$

and

$$\lambda_2(\rho_*, u_* + \int_0^t \alpha(s)\,ds) = \lambda_2(\rho_+, u_+ + \int_0^t \alpha(s)\,ds), \quad \frac{dx_2(t)}{dt} > \lambda_1(\rho_*, u_* + \int_0^t \alpha(s)\,ds)$$

for $t \geq 0$. As in [20], we use the Lax shock condition as an entropy criterion.

## 2.2. Delta-shock solution

From the previous analysis, we can see that the classical solution does not exist for any initial Riemann data. We denote the common boundary between regions $IV$ and $V$ by $S_\delta$. This notation stems from the fact that for right states $(\rho_+, u_+)$ in region $V$, we cannot connect $(\rho_-, u_-)$ and $(\rho_+, u_+)$ using classical shock waves $J_1$ and $J_2$. Therefore, it is necessary to employ delta-shock waves for right states in the region $V$ to solve the Riemann problem for system (1) with a source term.

We notice that, the curve $S_\delta$ has the line $u = u_-$ as its asymptotic line. When $(\rho_+, u_+) \in V$, we have $u_+ + \frac{\mu}{\rho_+} \leq u_-$. The characteristic curves from the initial data overlap in the domain

$$\Omega = \{(x,t) \; : \; \left(u_+ + \frac{\mu}{\rho_+}\right)t + \int_0^t \int_0^r \alpha(s)ds\,dr = x_2(t) \leq x \leq x_1(t) = u_-t + \int_0^t \int_0^r \alpha(s)ds\,dr, \; t > 0\}.$$

Therefore, singularities must occur in $\Omega$. Suppose $(\rho_+, u_+) \in V$. Let $(\rho_*, u_*)$ be any point on $J_1(\rho_-, u_-)$. For the point $(\rho_\bullet, u^\bullet)$ on $J_2(\rho_*, u_*)$, we have

$$u^\bullet + \frac{\mu}{\rho_+} + \int_0^t \alpha(s)ds = u_* + \frac{\mu}{\rho_*} + \int_0^t \alpha(s)ds = u_- + \frac{\mu}{\rho_*} + \int_0^t \alpha(s)ds$$

$$\geq u_+ + \frac{\mu}{\rho_+} + \frac{\mu}{\rho_*} + \int_0^t \alpha(s)ds.$$

Thus, $u^\bullet + \int_0^t \alpha(s)ds \geq u_+ + \frac{\mu}{\rho_+} + \int_0^t \alpha(s)ds$ for $\rho_* > 0$ and $(\rho_+, u_+)$ cannot lie on any curve in F, the singularity is impossible to be a jump with finite amplitude because the Rankine-Hugoniot conditions is not satisfied on the bounded jump. When $\alpha(\cdot) \equiv const.$, the solution for the Riemann problem (1) is a type of nonlinear hyperbolic waves called *delta-shock wave* (see [1]). Therefore, when $u_+ + \frac{\mu}{\rho_+} \leq u_-$, we suggest that the solution to the Riemann problem is also a delta-shock type wave defined by the speed $\sigma$ satisfying

$$u_+ + \frac{\mu}{\rho_+} + \int_0^t \alpha(s)ds \leq \sigma(t) \leq u_- + \int_0^t \alpha(s)ds, \qquad \text{for all } t \geq 0. \tag{13}$$

Next, we recall the following

**Definition 1.** Given $a, b \in \mathbb{R}$, $a < b$, and $w \in L^1((a,b))$, a two-dimensional weighted delta function $w(\cdot)\,\delta_L$ supported on a smooth curve $L = \{(x(s), t(s)) \; : \; a < s < b\}$, is defined as follows

$$< w(\cdot)\delta_L, \phi(\cdot, \cdot) >:= \int_a^b w(s)\phi(x(s), t(s))ds, \quad \text{for any } \phi \in C_0^\infty(\mathbb{R} \times [0, \infty)).$$

Now, we define a delta-shock wave solution for the Keyfitz-Kranzer system (1) with initial data (3).

**Definition 2.** A distribution pair $(\rho, u)$ is called a *delta-shock wave solution* of (1) and (3), in the sense of distributions, if there exists a smooth curve $L$ and a function $w \in C^1(L)$ such that $\rho$ and $u$ are represented as

$$\rho = \widehat{\rho}(x,t) + w\delta_L \quad \text{and} \quad u = \widehat{u}(x,t),$$

$\widehat{\rho}, \widehat{u} \in L^\infty(\mathbb{R} \times (0, \infty); \mathbb{R})$ and

$$\begin{cases} \langle \rho, \phi_t \rangle + \langle \rho u, \phi_x \rangle = 0, \\ \langle \rho u, \phi_t \rangle + \langle \rho u^2, \phi_x \rangle + \int_0^\infty \int_\mathbb{R} \mu u \phi_x \, dx dt = -\langle \alpha \rho, \phi \rangle \end{cases} \tag{14}$$

for all test functions $\phi \in C_0^\infty(\mathbb{R} \times (0, \infty))$, where $u|_L = u_\delta(t)$, and

$$\langle \rho, \phi \rangle := \int_0^\infty \int_\mathbb{R} \widehat{\rho}\phi \, dx dt + < w\delta_L, \phi >,$$

$$\langle \rho G(u), \phi \rangle := \int_0^\infty \int_\mathbb{R} \widehat{\rho}G(\widehat{u})\phi \, dx dt + < wG(u_\delta)\delta_L, \phi >,$$

$$\langle \alpha \rho, \phi \rangle := \int_0^\infty \int_\mathbb{R} \alpha(t)\widehat{\rho}\phi \, dx dt + < \alpha(t)w\delta_L, \phi > .$$

Then, with the previous definitions we will find a solution with discontinuity $x = x_\delta(t)$ for the system (1) of the form

$$(\rho(x,t), u(x,t)) = \begin{cases} (\rho_-(x,t), u_-(x,t)), & \text{if } x < x_\delta(t), \\ (w(t)\delta_L, u_\delta(t)), & \text{if } x = x_\delta(t), \\ (\rho_+(x,t), u_+(x,t)), & \text{if } x > x_\delta(t), \end{cases} \tag{15}$$

where $\rho_\pm(x,t)$, $u_\pm(x,t)$ are piecewise smooth solutions of system (1), $\delta(\cdot)$ is the Dirac measure supported on the curve $x_\delta \in C^1$, and the functions $x_\delta(t)$, $w(t)$, and $u_\delta(t)$ are to be determined. On the other hand, we define

$$\frac{1}{\rho} := \begin{cases} \frac{1}{\rho_-}, & \text{if } x < x_\delta(t), \\ 0, & \text{if } x = x_\delta(t), \\ \frac{1}{\rho_+}, & \text{if } x > x_\delta(t). \end{cases}$$

Now, we have the following

**Theorem 3.** *If the curves $x_\delta(t)$, $w(t)$, and $u_\delta(t)$ solve*

$$\begin{cases} \frac{dx_\delta(t)}{dt} = u_\delta(t), \\ \frac{dw(t)}{dt} = -u_\delta(t)[\![\rho]\!] + [\![\rho u]\!], \\ \frac{d(w(t)u_\delta(t))}{dt} = -u_\delta(t)[\![\rho u]\!] + [\![\rho u \left(u + \frac{\mu}{\rho}\right)]\!] + \alpha(t)w(t), \end{cases} \tag{16}$$

*then the solution $(\rho(x,t), u(x,t))$ defined in (15) satisfies (1) in the sense of distributions.*

Here $[\![q]\!] := q(x_\delta(t)-, t) - q(x_\delta(t)+, t)$.

PROOF. If equation (16) holds, then for any test functions $\phi \in C_0^\infty(\mathbb{R} \times (0, \infty))$, we obtain

$$\langle \rho u, \phi_t \rangle + \langle \rho u^2, \phi_x \rangle + \int_0^\infty \int_\mathbb{R} \mu u \phi_x \, dx dt$$

$$= \int_0^\infty \int_\mathbb{R} \left(\rho u \phi_t + \rho u \left(u + \frac{\mu}{\rho}\right)\phi_x\right) dx dt + \int_0^\infty (w(t)u_\delta(t)\phi_t + w(t)u_\delta(t)u_\delta(t)\phi_x)dt$$

$$= \int_0^\infty \int_{-\infty}^{x_\delta(t)} \left(\rho_- u_- \phi_t + \rho_- u_- \left(u_- + \frac{\mu}{\rho_-}\right)\phi_x\right) dx dt$$

$$+ \int_0^\infty \int_{x_\delta(t)}^\infty \left(\rho_+ u_+ \phi_t + \rho_+ u_+ \left(u_+ + \frac{\mu}{\rho_+}\right)\phi_x\right) dx dt + \int_0^\infty w(t)u_\delta(t)(\phi_t + u_\delta(t)\phi_x)dt$$

$$= -\oint -\left(\rho_- u_- \left(u_- + \frac{\mu}{\rho_-}\right)\phi\right)dt + (\rho_- u_- \phi)dx + \oint -\left(\rho_+ u_+ \left(u_+ + \frac{\mu}{\rho_+}\right)\phi\right)dt + (\rho_+ u_+ \phi)dx$$

$$- \int_0^\infty \int_\mathbb{R} \alpha(t)\rho\phi \, dx dt - \int_0^\infty \phi \frac{d(w(t)u_\delta(t))}{dt}dt$$

$$= \int_0^\infty \phi \left([\![\rho u \left(u + \frac{\mu}{\rho}\right)]\!] - u_\delta(t)[\![\rho u]\!] - \frac{d(w(t)u_\delta(t))}{dt}\right)dt \int_0^\infty \int_\mathbb{R} \alpha(t)\rho\phi \, dx dt - \int_0^\infty \alpha(t)w(t)\phi dt$$

$$+ \int_0^\infty \alpha(t)w(t)\phi dt$$

$$= \int_0^\infty \phi \left(-u_\delta(t)[\![\rho u]\!] + [\![\rho u \left(u + \frac{\mu}{\rho}\right)]\!] - \frac{d(w(t)u_\delta(t))}{dt} + \alpha(t)w(t)\right)dt \langle \alpha\rho, \phi \rangle = -\langle \alpha\rho, \phi \rangle$$

which implies the second equation of (17). A completely similar argument leads to the first equation of (17).

The relations given by (16) are called the *generalized Rankine-Hugoniot conditions* for the system (1).

We also define a delta-shockwave solution for the system (6).

**Definition 3.** A distribution pair $(\rho, u)$ is called a *delta-shock wave solution* of (6), in the sense of distributions, if there exists a smooth curve $L$ and a function $w \in C^1(L)$, such that, $\rho$ and $u$ are represented as

$$\rho = \widehat{\rho}(x,t) + w\delta_L \quad \text{and} \quad \widetilde{u} = \widehat{u}(x,t),$$

$\widehat{\rho}, \widehat{u} \in L^\infty(\mathbb{R} \times (0, \infty); \mathbb{R})$ and

$$\begin{cases} \langle \rho, \phi_t \rangle + \langle \rho(\widetilde{u} + \int_0^t \alpha(s)ds), \phi_x \rangle = 0, \\ \langle \rho\widetilde{u}, \phi_t \rangle + \langle \rho\widetilde{u}(\widetilde{u} + \int_0^t \alpha(s)ds), \phi_x \rangle + \int_0^\infty \int_\mathbb{R} \mu\widetilde{u}\phi_x \, dx dt = 0, \end{cases} \tag{17}$$

for all test functions $\phi \in C_0^\infty(\mathbb{R} \times (0, \infty))$, where $u|_L = u_\delta(t)$, and

$$\langle \rho, \phi \rangle := \int_0^\infty \int_\mathbb{R} \widehat{\rho}\phi \, dx dt + \langle w\delta_L, \phi \rangle,$$

$$\langle \rho G(\widetilde{u}), \phi \rangle := \int_0^\infty \int_\mathbb{R} \widehat{\rho}G(\widehat{u})\phi \, dx dt + \langle wG(u_\delta)\delta_L, \phi \rangle.$$

Thus, we need to find a solution for problem (6)-(7) of the form

$$(\rho(x,t), \widetilde{u}(x,t)) = \begin{cases} (\rho_-, u_-), & \text{if } x < x_\delta(t), \\ (w(t)\delta_L, \widetilde{u}_\delta(t)), & \text{if } x = x_\delta(t), \\ (\rho_+, u_+), & \text{if } x > x_\delta(t), \end{cases} \tag{18}$$

where $\rho_\pm$, $u_\pm$ are constants, $\delta(\cdot)$ is the Dirac measure supported on the curve $x_\delta(t) \in C^1$, and $x_\delta(t)$, $w(t)$ and $\widetilde{u}_\delta(t)$ are to be determined. Then, we have the following result

**Lemma 4.** *If the curves $x_\delta(t)$, $w(t)$, and $\widetilde{u}_\delta(t)$ solve*

$$\begin{cases} \frac{dx_\delta(t)}{dt} = \widetilde{u}_\delta(t) + \int_0^t \alpha(s)ds, \\ \frac{dw(t)}{dt} = -(\widetilde{u}_\delta(t) + \int_0^t \alpha(s)ds)[\rho] + [\rho(\widetilde{u} + \int_0^t \alpha(s)ds)], \\ \frac{d(w(t)\widetilde{u}_\delta(t))}{dt} = -(\widetilde{u}_\delta(t) + \int_0^t \alpha(s)ds)[\rho\widetilde{u}] + [\rho\widetilde{u}(\widetilde{u} + \int_0^t \alpha(s)ds + \frac{\mu}{\rho})], \end{cases} \tag{19}$$

*then the solution $(\rho(x,t), \widetilde{u}(x,t))$ defined in (18) satisfies the problem (6)-(7) in the sense of distributions.*

The proof is similar to that one obtained in Theorem 3, hence we omit it. From the previous lemma, the Riemann problem to (6) is reduced to find $x_\delta(t)$, $w(t)$, and $\widetilde{u}_\delta(t)$, such that, they solve the system (19) with the initial data (at $t = 0$), that is

$$x_\delta(0) = 0, \qquad w(0) = 0, \quad \text{and} \quad \widetilde{u}_\delta(0) = \eta = const..$$

Under the following entropy condition

$$\widetilde{\lambda}_2(\rho_+, u_+) = u_+ + \frac{\mu}{\rho_+} + \int_0^t \alpha(s)ds \leq \sigma(t) = \frac{dx_\delta(t)}{dt} \leq u_- + \int_0^t \alpha(s)ds = \widetilde{\lambda}_1(\rho_-, u_-), \quad \forall t \geq 0, \tag{20}$$

from the second equation of (19), and using the inequality (20) we obtain

$$\frac{dw(t)}{dt} = \rho_-\left(u_- - \widetilde{u}_\delta(t)\right) + \rho_+(\widetilde{u}_\delta(t) - u_+) > 0.$$

In addition, by integrating 0 to $t$ the second equation of (19) we get

$$w(t) = \rho_- \left( u_- t + \int_0^t \int_0^r \alpha(s) \, ds \, dr - x_\delta(t) \right) + \rho_+ \left( x_\delta(t) - u_+ t + \int_0^t \int_0^r \alpha(s) \, ds \, dr \right) \geq 0.$$

Moreover, by integrating 0 to $t$ each equation in (19) the following system results

$$\begin{cases} x_\delta(t) = \int_0^t \widetilde{u}_\delta(r) \, dr + \int_0^t \int_0^r \alpha(s) \, ds \, dr, \\ w(t) = -[\rho] \int_0^t \widetilde{u}_\delta(r) \, dr + [\rho \widetilde{u}] t, \\ w(t) \widetilde{u}_\delta(t) = -[\rho \widetilde{u}] \int_0^t \widetilde{u}_\delta(r) \, dr + [\rho \widetilde{u}^2 + \mu \widetilde{u}] t. \end{cases} \tag{21}$$

From the second and third equations of system (21), we have

$$[\rho] \widetilde{u}_\delta(t) \int_0^t \widetilde{u}_\delta(r) \, dr - [\rho \widetilde{u}] \widetilde{u}_\delta(t) t - [\rho \widetilde{u}] \int_0^t \widetilde{u}_\delta(r) \, dr + [\rho \widetilde{u}^2 + \mu \widetilde{u}] t = 0,$$

or

$$\frac{1}{2} \frac{d}{dt} \left( [\rho] \left( \int_0^t \widetilde{u}_\delta(r) \, dr \right)^2 - 2[\rho \widetilde{u}] t \int_0^t \widetilde{u}_\delta(r) \, dr + [\rho \widetilde{u}^2 + \mu \widetilde{u}] t^2 \right) = 0 \tag{22}$$

and integrating (22) from 0 to $t$, we obtain the following equality

$$[\rho] \left( \int_0^t \widetilde{u}_\delta(r) \, dr \right)^2 - 2[\rho \widetilde{u}] t \int_0^t \widetilde{u}_\delta(r) \, dr + [\rho \widetilde{u}^2 + \mu \widetilde{u}] t^2 = 0. \tag{23}$$

From (23), one can find that $\widetilde{u}_\delta(t)$ is a constant, we say that $\widetilde{u}_\delta(t) = u_\delta = const.$, and $\int_0^t \widetilde{u}_\delta(r) \, dr = u_\delta t$. Then, we can rewrite (23) as $[\rho] u_\delta^2 - 2[\rho \widetilde{u}] u_\delta + [\rho \widetilde{u}^2 + \mu \widetilde{u}] = 0$. Therefore, if $[\rho] = \rho_- - \rho_+ = 0$, then we have

$$u_\delta = \frac{1}{2} \left( u_- + u_+ + \frac{\mu}{\rho_+} \right), \quad x_\delta(t) = \frac{1}{2} \left( u_- + u_+ + \frac{\mu}{\rho_+} \right) t + \int_0^t \int_0^r \alpha(s) \, ds \, dr,$$

$$w(t) = \rho_+ (u_- - u_+) t, \tag{24}$$

which satisfies the entropy condition (20). When $[\rho] = \rho_- - \rho_+ \neq 0$, the discriminant of the quadratic equation (23) is

$$\Delta = 4([\rho \widetilde{u}]^2 - [\rho][\rho \widetilde{u}^2 + \mu \widetilde{u}]) = 4\rho_- \rho_+ (u_- - u_+) \left( u_- + \frac{\mu}{\rho_-} - u_- - \frac{\mu}{\rho_+} \right) > 0,$$

and then we can find

$$\begin{cases} u_\delta = \frac{[\rho \widetilde{u}] - \sqrt{[\rho \widetilde{u}]^2 - [\rho][\rho \widetilde{u}^2 + \mu \widetilde{u}]}}{[\rho]}, \\ x_\delta(t) = \frac{[\rho \widetilde{u}] - \sqrt{[\rho \widetilde{u}]^2 - [\rho][\rho \widetilde{u}^2 + \mu \widetilde{u}]}}{[\rho]} t + \int_0^t \int_0^r \alpha(s) \, ds \, dr, \\ w(t) = \sqrt{[\rho \widetilde{u}]^2 - [\rho][\rho \widetilde{u}^2 + \mu \widetilde{u}]} \, t, \end{cases} \tag{25}$$

or

$$\begin{cases} u_\delta = \frac{[\rho \widetilde{u}] + \sqrt{[\rho \widetilde{u}]^2 - [\rho][\rho \widetilde{u}^2 + \mu \widetilde{u}]}}{[\rho]}, \\ x_\delta(t) = \frac{[\rho \widetilde{u}] + \sqrt{[\rho \widetilde{u}]^2 - [\rho][\rho \widetilde{u}^2 + \mu \widetilde{u}]}}{[\rho]} t + \int_0^t \int_0^r \alpha(s) \, ds \, dr, \\ w(t) = -\sqrt{[\rho \widetilde{u}]^2 - [\rho][\rho \widetilde{u}^2 + \mu \widetilde{u}]} \, t. \end{cases} \tag{26}$$

Next, using the entropy condition (20), we choose the admissible solution between (25) and (26). We observe that, for the solution given in (25) we have

$$\frac{dx_\delta(t)}{dt} - \widetilde{\lambda}_1(\rho_-, \widetilde{u}_-) = \frac{[\rho\widetilde{u}] - \sqrt{[\rho\widetilde{u}]^2 - [\rho][\rho\widetilde{u}^2 + \mu\widetilde{u}]}}{[\rho]} + \int_0^t \alpha(s)ds - u_- - \int_0^t \alpha(s)ds$$

$$= \frac{[\rho\widetilde{u}] - \sqrt{[\rho\widetilde{u}]^2 - [\rho][\rho\widetilde{u}^2 + \mu\widetilde{u}]} - u_-[\rho]}{[\rho]} = \frac{\rho_+[u] - \sqrt{[\rho\widetilde{u}]^2 - [\rho][\rho\widetilde{u}^2 + \mu\widetilde{u}]}}{[\rho]}$$

$$= \frac{-\rho_+[u]^2 + \mu[u]}{\rho_+[u] + \sqrt{[\rho\widetilde{u}]^2 - [\rho][\rho\widetilde{u}^2 + \mu\widetilde{u}]}} = \frac{-\rho_+[u]([u] - \frac{\mu}{\rho_+})}{\rho_+[u] + \sqrt{[\rho\widetilde{u}]^2 - [\rho][\rho\widetilde{u}^2 + \mu\widetilde{u}]}} \leq 0$$

and

$$\frac{dx_\delta(t)}{dt} - \widetilde{\lambda}_2(\rho_-, \widetilde{u}_-) = \frac{\rho_-([\widetilde{u}] - \mu/\rho_+)([\widetilde{u}] + \mu/\rho_- - \mu/\rho_+)}{\rho_-([\widetilde{u}] + \mu/\rho_- - \mu/\rho_+) + \sqrt{[\rho\widetilde{u}]^2 - [\rho][\rho\widetilde{u}^2 + \mu\widetilde{u}]}} \geq 0,$$

which implies that the entropy condition (20) is valid. Now, for the solution (26), if $\rho_- > \rho_+$ then

$$\frac{dx_\delta(t)}{dt} - \widetilde{\lambda}_1(\rho_-, \widetilde{u}_-) = \frac{\rho_+[\widetilde{u}] + \sqrt{[\rho\widetilde{u}]^2 - [\rho][\rho\widetilde{u}^2 + \mu\widetilde{u}]}}{[\rho]} > 0,$$

and when $\rho_- < \rho_+$, we obtain

$$\frac{dx_\delta(t)}{dt} - \widetilde{\lambda}_2(\rho_-, \widetilde{u}_-) = \frac{\rho_-([\widetilde{u}] + \mu/\rho_- - \mu/\rho_+) + \sqrt{[\rho\widetilde{u}]^2 - [\rho][\rho\widetilde{u}^2 + \mu\widetilde{u}]}}{[\rho]} < 0,$$

and the solution (26) does not satisfy the entropy condition (20).

From the previous analysis, we have proved the following result

**Proposition 5.** *Given left and right constant states $(\rho_-, u_-)$ and $(\rho_+, u_+)$, respectively, such that $u_+ + \mu/\rho_+ \leq u_-$, then the Riemann problem (6)-(7) admits a unique entropy solution in the sense of the Definition 3. This solution is of the form (18) where $\widetilde{u}_\delta = u_\delta$, $x_\delta(t)$, and $w(t)$ are shown in (24) if $[\rho] = 0$ or in (25) if $[\rho] \neq 0$.*

Now, we have the following result for the Riemann problem to nonsymmetric Keyfitz-Kranzer system (1).

**Theorem 6.** *Given left and right constant states $(\rho_-, u_-)$ and $(\rho_+, u_+)$, respectively, such that $u_+ + \mu/\rho_+ \leq u_-$, the Riemann problem for the nonsymmetric Keyfitz-Kranzer system with time-dependent Coulomb-like friction term (1) and initial data (3), admits a unique entropy solution in the sense of the Definition 2. This solution is of the form*

$$(\rho(x,t), u(x,t)) = \begin{cases} (\rho_-, u_- + \int_0^t \alpha(s)ds), & \text{if } x < x_\delta(t), \\ (w(t)\delta(x_\delta(t) - x), u_\delta(t)), & \text{if } x = x_\delta(t), \\ (\rho_+, u_+ + \int_0^t \alpha(s)ds), & \text{if } x > x_\delta(t), \end{cases} \tag{27}$$

*where $u_\delta(t)$, $x_\delta(t)$, and $w(t)$ are given by*

$$\begin{cases} u_\delta(t) = \frac{1}{2}\left(u_- + u_+ + \frac{\mu}{\rho_+}\right) + \int_0^t \alpha(s)ds, \\ x_\delta(t) = \frac{1}{2}\left(u_- + u_+ + \frac{\mu}{\rho_+}\right)t + \int_0^t \int_0^r \alpha(s)ds\,dr, \\ w(t) = \rho_+(u_- - u_+)t, \end{cases} \tag{28}$$

*when $\rho_- - \rho_+ = 0$, or given by*

$$\begin{cases} u_\delta(t) = \frac{[\rho u] - \sqrt{[\rho u]^2 - [\rho][\rho u^2 + \mu u]}}{[\rho]} + \int_0^t \alpha(s)\,ds, \\ x_\delta(t) = \frac{[\rho u] - \sqrt{[\rho u]^2 - [\rho][\rho u^2 + \mu u]}}{[\rho]} t + \int_0^t \int_0^r \alpha(s)\,ds\,dr, \\ w(t) = \sqrt{[\rho u]^2 - [\rho][\rho u^2 + \mu u]}\, t, \end{cases} \tag{29}$$

*when $\rho_- - \rho_+ \neq 0$.*

PROOF. In order to solve (16), we consider the initial data (at $t = 0$), $x_\delta(0) = 0$, $w(0) = 0$, and $u_\delta(0) = u_\delta$, under the following over-compressive entropy condition

$$u_+ + \frac{\mu}{\rho_+} + \int_0^t \alpha(s)\,ds \leq \frac{dx_\delta(t)}{dt} = u_\delta(t) \leq u_- + \int_0^t \alpha(s)\,ds, \qquad \text{for all } t \geq 0.$$

From (5), for the initial data given by the left and right constant states $(\rho_-, u_-)$ and $(\rho_+, u_+)$, respectively, we have

$$\begin{cases} [\![\rho]\!] = \rho_- - \rho_+ = [\rho], \\ [\![\rho u]\!] = \rho_-(u_- + \int_0^t \alpha(s)\,ds) - \rho_+(u_+ + \int_0^t \alpha(s)\,ds) = [\rho u] + [\rho]\int_0^t \alpha(s)\,ds, \\ [\![\rho u^2 + \mu u]\!] = [\rho u^2 + \mu u] + 2[\rho u]\int_0^t \alpha(s)\,ds + [\rho](\int_0^t \alpha(s)\,ds)^2, \end{cases}$$

and from (16) we obtain

$$\begin{cases} \frac{dx_\delta(t)}{dt} = u_\delta(t), \\ \frac{dw(t)}{dt} = -u_\delta(t)[\rho] + [\rho u] + [\rho]\int_0^t \alpha(s)\,ds, \\ \frac{d(w(t)u_\delta(t))}{dt} = -u_\delta(t)([\rho u] + [\rho]\int_0^t \alpha(s)\,ds) + [\rho u^2 + \mu u] + 2[\rho u]\int_0^t \alpha(s)\,ds + [\rho](\int_0^t \alpha(s)\,ds)^2 + \alpha(t)w(t). \end{cases} \tag{30}$$

From the second and third equations of the system (30), it follows that

$$\frac{d(w(t)(u_\delta(t) - \int_0^t \alpha(s)\,ds))}{dt} = (u_\delta(t) - \int_0^t \alpha(s)\,ds)[\rho u] + [\rho u^2 + \mu u],$$

and therefore we obtain

$$\begin{cases} x_\delta(t) = \int_0^t u_\delta(r)\,dr, \\ w(t) = -[\rho]\int_0^t (u_\delta(r) - \int_0^r \alpha(s)\,ds)\,dr + [\rho u]t, \\ w(t)(u_\delta(t) - \int_0^t \alpha(s)\,ds) = -[\rho u]\int_0^t (u_\delta(r) - \int_0^r \alpha(s)\,ds)\,dr + [\rho u^2 + \mu u]t. \end{cases}$$

which implies that

$$[\rho]\frac{1}{2}\frac{d}{dt}\left(\int_0^t (u_\delta(r) - \int_0^r \alpha(s)\,ds)\,dr\right)^2 - [\rho u]\frac{d}{dt}\left(t\int_0^t (u_\delta(r) - \int_0^r \alpha(s)\,ds)\,dr\right) + [\rho u^2 + \mu u]\frac{1}{2}\frac{d}{dt}(t^2) = 0$$

or

$$[\rho]\left(\int_0^t (u_\delta(r) - \int_0^r \alpha(s)\,ds)\,dr\right)^2 2[\rho u]t\int_0^t (u_\delta(r) - \int_0^r \alpha(s)\,ds)\,dr + +[\rho u^2 + \mu u]t^2 = 0. \tag{31}$$

From (31), we find that $u_\delta(t) - \int_0^t \alpha(s)\,ds$ is a constant, we say $u_\delta(t) - \int_0^t \alpha(s)\,ds = u_\delta = const.$ or $u_\delta(t) = u_\delta + \int_0^t \alpha(s)\,ds$. Similar as the results to obtain the Proposition 5, one can show that $u_\delta = \frac{1}{2}(u_- + u_+ + \mu/\rho_+)$ for $[\rho] = 0$ and $u_\delta = \frac{[\rho u] - \sqrt{[\rho u]^2 - [\rho][\rho u^2 + \mu u]}}{[\rho]}$ for $[\rho] \neq 0$, under the entropy (13). In addition, we have that $x_\delta(t) = u_\delta t + \int_0^t \int_0^r \alpha(s)\,ds\,dr$ and $w(t) = (-[\rho]u_\delta + [\rho u])t$.

Now, applying Theorem 3, we have Eqs. (28) and (29) are satisfied, and we conclude that the solution to the Riemann problem (1) is given by (27).

## 3. The Pressureless system with time-dependent Coulomb-like friction term

In this section, we consider a pressureless gas system with a time-dependent Coulomb-like friction term (2) with initial data (3). Applying a suitable transformation together with a similar variable, Riemann solutions are constructed. Indeed, two contact discontinuities are obtained for $u_- < u_+$, and if $u_- > u_+$, then a delta-shock solution is established.

To begin, using the so-called Dafermos regularization, we have the following viscous system associated with the system (2),

$$\begin{cases} \rho_t^\varepsilon + (\rho^\varepsilon u^\varepsilon)_x = 0, \\ (\rho^\varepsilon u^\varepsilon)_t + (\rho^\varepsilon (u^\varepsilon)^2)_x = \alpha(t)\rho^\varepsilon + \varepsilon t u_{xx}^\varepsilon, \end{cases} \tag{32}$$

and with the transformation

$$u^\varepsilon(x,t) = \widetilde{u}^\varepsilon(x,t) + \int_0^t \alpha(s)\,ds \tag{33}$$

the previous system can be written as

$$\begin{cases} \rho_t^\varepsilon + \left(\rho^\varepsilon \left(\widetilde{u}^\varepsilon + \int_0^t \alpha(s)\,ds\right)\right)_x = 0, \\ (\rho^\varepsilon \widetilde{u}^\varepsilon)_t + \left(\rho^\varepsilon \widetilde{u}^\varepsilon \left(\widetilde{u}^\varepsilon + \int_0^t \alpha(s)\,ds\right)\right)_x = \varepsilon t \widetilde{u}_{xx}^\varepsilon. \end{cases} \tag{34}$$

Now, introducing the following similar variable

$$\xi = \frac{x - \int_0^t \int_0^r \alpha(s)\,ds\,dr}{t} \tag{35}$$

the system (34) can be written as

$$\begin{cases} -\xi \rho_\xi^\varepsilon + (\rho^\varepsilon \widetilde{u}^\varepsilon)_\xi = 0, \\ -\xi (\rho^\varepsilon \widetilde{u}^\varepsilon)_\xi + (\rho^\varepsilon (\widetilde{u}^\varepsilon)^2)_\xi = \varepsilon \widetilde{u}_{\xi\xi}^\varepsilon, \end{cases} \tag{36}$$

and the initial data (7) changes to the boundary condition

$$(\rho(\pm\infty), \widetilde{u}(\pm\infty)) = (\pm\rho, \pm u). \tag{37}$$

According to Section 2.2 in [29], there exists a solution to the boundary problem (36)-(37). In addition, the limit solutions as $\varepsilon \to 0$ are given in the next two propositions, and their proofs are given in Section 2.3 in [29].

**Proposition 7.** *Suppose $u_- < u_+$ and let $(\rho, \widetilde{u})$ be the solution to the problem* (36)-(37). *Then,*

$$\lim_{\varepsilon \to 0}(\rho^\varepsilon(\xi), \widetilde{u}^\varepsilon(\xi)) = \begin{cases} (\rho_-, u_-), & \text{if } \xi < u_-, \\ (0, \xi), & \text{if } u_- < \xi < u_+, \\ (\rho_+, u_+), & \text{if } \xi > u_+. \end{cases}$$

**Proposition 8.** *Suppose $u_- > u_+$ and let $(\rho, \widetilde{u})$ be the solution to the problem* (36)-(37). *Then,*

$$\lim_{\varepsilon \to 0}(\rho^\varepsilon(\xi), \widetilde{u}^\varepsilon(\xi)) = \begin{cases} (\rho_-, u_-), & \text{if } \xi < u_\delta, \\ (w_\delta \delta(\xi - u_\delta), u_\delta), & \text{if } \xi = u_\delta, \\ (\rho_+, u_+), & \text{if } \xi > u_\delta, \end{cases}$$

*where*

$$w_\delta = \sqrt{\rho_- \rho_+}(u_- - u_+), \quad u_\delta = \frac{\sqrt{\rho_-}\, u_- + \sqrt{\rho_+}\, u_+}{\sqrt{\rho_-} + \sqrt{\rho_+}}.$$

*Moreover, $\rho^\varepsilon$ converges in the sense of distribution to the sum of a step function and a Dirac measure $\delta$ with weight $w_\delta = -u_\delta(\rho_- - \rho_+) + (\rho_- u_- - \rho_+ u_+)$ and $\rho^\varepsilon u^\varepsilon$ converges in the sense of distribution to the sum of a step function and a Dirac measure $\delta$ with weight $-u_\delta(\rho_- u_- - \rho_+ u_+) + (\rho_- u_-^2 - \rho_+ u_+^2)$.*

One observes that, as $\varepsilon \to 0$, (33) converges to (5). Thus, using the similar variable (35), from Propositions 7 and (8) we have the following results for the pressureless system with a time-dependent Coulomb-like friction term (2).

**Theorem 9.** *Suppose $u_- < u_+$ and let $(\rho^\varepsilon, u^\varepsilon)$ be the solution to the problem* (32)-(3). *Then, the pressureless system with a time-dependent Coulomb-like friction term* (2) *with initial data* (3) *has the following solution*

$$(\rho(x,t), u(x,t)) = \lim_{\varepsilon \to 0}(\rho^\varepsilon(x,t), u^\varepsilon(x,t)) = \begin{cases} (\rho_-, u_- + \int_0^t \alpha(s)\,ds), & \text{if } x < x_-(t), \\ (0, \frac{x - \int_0^t \int_0^r \alpha(s)\,ds\,dr}{t} + \int_0^t \alpha(s)\,ds), & \text{if } x_-(t) < x < x_+(t), \\ (\rho_+, u_+ + \int_0^t \alpha(s)\,ds), & \text{if } x > x_+(t), \end{cases}$$

*where $x_\pm(t) = u_\pm t + \int_0^t \int_0^r \alpha(s)\,ds\,dr$.*

**Theorem 10.** *Suppose $u_- > u_+$ and let $(\rho^\varepsilon, u^\varepsilon)$ be the solution to the problem* (32)-(3). *Then, the pressureless system with a time-dependent Coulomb-like friction term* (2) *with initial data* (3) *has the following entropy solution*

$$(\rho(x,t), u(x,t)) = \lim_{\varepsilon \to 0}(\rho^\varepsilon(x,t), u^\varepsilon(x,t)) = \begin{cases} (\rho_-, u_- + \int_0^t \alpha(s)\,ds), & \text{if } x < x_\delta(t), \\ (w(t)\delta(x_\delta(t) - x), u_\delta(t)), & \text{if } x = x_\delta(t), \\ (\rho_+, u_+ + \int_0^t \alpha(s)\,ds), & \text{if } x > x_\delta(t), \end{cases}$$

*where*

$$w(t) = \sqrt{\rho_- \rho_+}(u_- - u_+)t, \; x_\delta(t) = \frac{\sqrt{\rho_-}\, u_- + \sqrt{\rho_+}\, u_+}{\sqrt{\rho_-} + \sqrt{\rho_+}}t + \int_0^t \int_0^r \alpha(s)\,ds\,dr,$$

*and $u_\delta(t) = \frac{\sqrt{\rho_-}\, u_- + \sqrt{\rho_+}\, u_+}{\sqrt{\rho_-} + \sqrt{\rho_+}} + \int_0^t \alpha(s)\,ds$.*

## 4. The vanishing pressure limit of Riemann solutions to System (1)

In this section, we consider the vanishing pressure limit, that is $\mu \to 0$, of the Riemann solutions to the problem (1), (3). According to the relations between $u_-$ and $u_+$, we divide our discussion into the following two cases: $u_- < u_+$ and $u_- > u_+$.

### 4.1. Similar variable and Riemann problem to the nonsymmetric Keyfitz-Kranzer System

Using the similar variable (35), the system (6) can be written as

$$\begin{cases} -\xi \rho_\xi + (\rho\widetilde{u})_\xi = 0, \\ -\xi(\rho\widetilde{u})_\xi + (\rho\widetilde{u}^2 + \mu\widetilde{u})_\xi = 0, \end{cases} \tag{38}$$

and the initial data (7) changes to the boundary condition

$$(\rho(\pm\infty), \widetilde{u}(\pm\infty)) = (\pm\rho, \pm u). \tag{39}$$

For $u_- < u_+ + \mu/\rho_+$, it is easy to solve the boundary problem (38)-(39) and using the similar variable (35) we have the explicit solution given in Proposition 1. This solution can be used together with the transformation (5) to obtain the classical solution given in Theorem 2. For the case $u_- > u_+ + \mu/\rho_+$, we consider the following viscosity system:

$$\begin{cases} \rho_t^\varepsilon + (\rho^\varepsilon u^\varepsilon)_x = 0, \\ (\rho^\varepsilon u^\varepsilon)_t + (\rho^\varepsilon (u^\varepsilon)^2 + \mu u^\varepsilon)_x = \alpha(t)\rho^\varepsilon + \varepsilon t u_{xx}^\varepsilon, \end{cases}$$

and with the transformation (33), the previous system can be written as

$$\begin{cases} \rho_t^\varepsilon + (\rho^\varepsilon(\widetilde{u}^\varepsilon + \int_0^t \alpha(s)ds))_x = 0, \\ (\rho^\varepsilon \widetilde{u}^\varepsilon)_t + \left(\rho^\varepsilon \widetilde{u}^\varepsilon(\widetilde{u}^\varepsilon + \int_0^t \alpha(s)ds) + \mu\widetilde{u}^\varepsilon\right)_x = \varepsilon t \widetilde{u}_{xx}^\varepsilon. \end{cases} \tag{40}$$

Using the similar variable (35), the system (40) is rewritten as

$$\begin{cases} -\xi\rho_\xi^\varepsilon + (\rho^\varepsilon \widetilde{u}^\varepsilon)_\xi = 0, \\ -\xi(\rho^\varepsilon \widetilde{u}^\varepsilon)_\xi + (\rho^\varepsilon(\widetilde{u}^\varepsilon)^2 + \mu\widetilde{u}^\varepsilon)_\xi = \varepsilon\widetilde{u}_{\xi\xi}^\varepsilon \end{cases} \tag{41}$$

and the initial data (7) changes to the boundary condition

$$(\rho(\pm\infty), \widetilde{u}(\pm\infty)) = (\pm\rho, \pm u). \tag{42}$$

The existence of solutions to the system (41) with boundary condition (42) was shown in [10]. Moreover, for $u_- > u_+ + \mu/\rho_+$, the solution to the boundary problem (41)-(42) is given in the following

**Theorem 11.** *Suppose $u_- > u_+ + \mu/\rho_+$ and let $(\rho^\varepsilon(\xi), u^\varepsilon(\xi))$ be the solution to the problem* (41)-(42). *Then,*

$$\lim_{\varepsilon\to 0}(\rho^\varepsilon(\xi), \widetilde{u}^\varepsilon(\xi)) = \begin{cases} (\rho_-, u_-), & \text{if } \xi < u_\delta, \\ (w_0\delta, u_\delta), & \text{if } \xi = u_\delta, \\ (\rho_+, u_+), & \text{if } \xi > u_\delta, \end{cases}$$

*where $w_0$ and $u_\delta$ are determined uniquely by the entropy condition $u_+ + \mu/\rho_+ < u_\delta < u_-$ and*

$$\begin{cases} w_0 = -u_\delta(\rho_- - \rho_+) + (\rho_- u_- - \rho_+ u_+), \\ w_0 u_\delta = -u_\delta(\rho_- u_- - \rho_+ u_+) + (\rho_- u_-^2 - \rho_+ u_+^2 + \mu u_- - \mu u_+). \end{cases}$$

The proof of Theorem 11 can be found in [10] (see proof of Theorem 5 of that reference). Again, we use the similarity variable (35) with the transformation (5) and Theorem 11 to obtain the delta-shock solution given in Theorem 6.

## 4.2. Vanshing pressure limit of the Riemann solution of the Keyfitz-Kranzer System

In this section, we analyze the behavior when the Chaplygin pressure $-\mu/\rho$ approaches 0 and the convergence of the Riemann solutions of the problem (1)-(3). Thus, to avoid confusion in this section, we denote by $(\rho^\mu, u^\mu)$ the solutions for the system (1). We divide the analysis into two cases: $u_- < u_+$ and $u_- > u_+$.

**Case 1** ($u_- < u_+$): In this case, $(\rho_+, u_+)$ belongs to $I \cup II$ in the phase plane. Thus, the Riemann solution to the problem (1)-(3) is given by (12) with $u_*^\mu = u_-$ and $\rho_*^\mu = \frac{\mu\rho_+}{\mu + \rho_+(u_+ - u_-)}$. Therefore, we have

$$\lim_{\mu\to 0}\rho_*^\mu = \lim_{\mu\to 0}\left(\frac{\mu\rho_+}{\mu + \rho_+(u_+ - u_-)}\right) = 0,$$

which indicates the occurrence of vacuum states. Moreover, from (10), the contact discontinuities are

$$x_1(t) = u_- t + \int_0^t \int_0^r \alpha(s)\,ds\,dr, \quad x_2(t) = \left(u_+ + \frac{\mu}{\rho_+}\right)t + \int_0^t \int_0^r \alpha(s)\,ds\,dr$$

which converge, respectively, to $x_-(t) = u_- t + \int_0^t \int_0^r \alpha(s)\,ds\,dr$ and $x_+(t) = u_+ t + \int_0^t \int_0^r \alpha(s)\,ds\,dr$ as $\mu \to 0$. In addition, as $u_- \le u_*^\mu \le u_+$, we have that $x_1 < x < x_2$ where $x = x(t) = u_*^\mu t + \int_0^t \int_0^r \alpha(s)\,ds\,dr$ and thus

$$\lim_{\mu \to 0} u_*^\mu = \lim_{\mu \to 0} \left(u_+ + \frac{\mu}{\rho_+} - \frac{\mu}{\rho_*}\right) = \frac{x}{t} - \frac{1}{t} \int_0^t \int_0^r \alpha(s)\,ds\,dr.$$

Therefore, the Riemann solution to the problem (1)-(3) converge to

$$\lim_{\mu \to 0}(\rho^\mu(x,t), u^\mu(x,t)) = \begin{cases} (\rho_-, u_- + \int_0^t \alpha(s)\,ds), & \text{if } x < x_-(t), \\ (0, \frac{x}{t} - \frac{1}{t} \int_0^t \int_0^r \alpha(s)\,ds\,dr), & \text{if } x_-(t) < x < x_+(t), \\ (\rho_+, u_+ + \int_0^t \alpha(s)\,ds), & \text{if } x > x_+(t), \end{cases}$$

where $x_\pm(t) = u_\pm t + \int_0^t \int_0^r \alpha(s)\,ds\,dr$.

**Case 2** ($u_- > u_+$): In this case, $(\rho_+, u_+)$ belongs to regions $III$, $IV$ or $V$. However, to analyze the behavior as $\mu \to 0$, we only consider cases in which $(\rho_+, u_+)$ belongs to regions $IV$ or $V$. Indeed, in the next lemma, we obtain that there is a constant $\mu_0 > 0$, such that, $\mu > \mu_0$ if $(\rho_+, u_+) \in IV$ and $0 < \mu < \mu_0$ if $(\rho_+, u_+) \in V$. Therefore, $\mu$ can tend to $\mu_0 > 0$ if $(\rho_+, u_+) \in IV$ (and a similar behavior occurs for the case when $(\rho_+, u_+) \in III$, hence it is not necessary to study it). The main goal in this section is to analyze the behavior as $\mu \to 0$, which is possible when $(\rho_+, u_+) \in V$.

**Lemma 12.** *Suppose $u_- > u_+$. Then, there is $\mu_0 > 0$ such that $(\rho_+, u_+) \in IV$ as $\mu > \mu_0$ and $(\rho_+, u_+) \in V$ as $\mu < \mu_0$. More explicitly, $\mu_0 = \rho_+(u_- - u_+)$.*

PROOF. If $(\rho_+, u_+) \in IV$, then $u_+ < u_- < u_+ + \mu/\rho_+$ and therefore $\mu_0 = \rho_+(u_- - u_+) < \mu$. If $(\rho_+, u_+) \in V$, $u_+ < u_-$ and $\mu/(u_- - u_+) < \rho_+$ or $\mu < \rho_+(u_- - u_+) = \mu_0$.

From the above lemma, we have that $(\rho_+, u_+) \in IV$ when $\mu_0 < \mu$. Therefore, from Theorem 2 we have $u_*^\mu = u_-$ and $\rho_*^\mu = \frac{\mu \rho_+}{\mu + \rho_+(u_+ - u_-)}$. Then, we obtain

$$\lim_{\mu \to \mu_0} \rho_*^\mu = \lim_{\mu \to \mu_0} \frac{\mu \rho_+}{\mu + \rho_+(u_+ - u_-)} = \lim_{\mu \to \mu_0} \frac{\mu \rho_+}{\mu - \mu_0} = \infty. \tag{43}$$

In the next lemma, we show that when $\mu \to \mu_0$, $\rho$ has the same singularity as a weighted Dirac delta function. Moreover, the curves $x_1$ and $x_2$ converge to a new type of nonlinear wave denoted by $\sigma$.

**Lemma 13.** *Suppose $u_- > u_+$ and $\mu > \mu_0 = \rho_+(u_- - u_+)$. Then,*

$$\lim_{\mu \to \mu_0} (u_*^\mu + \int_0^t \alpha(s)\,ds) = \lim_{\mu \to \mu_0} \frac{dx_1(t)}{dt} = \lim_{\mu \to \mu_0} \frac{dx_2(t)}{dt} = u_- + \int_0^t \alpha(s)\,ds =: \sigma(t), \tag{44}$$

$$\lim_{\mu \to \mu_0} \int_{x_1(t)}^{x_2(t)} \rho_*^\mu\,dx = \mu_0 t, \tag{45}$$

$$\lim_{\mu \to \mu_0} \int_{x_1(t)}^{x_2(t)} \rho_*^\mu (u_*^\mu + \int_0^t \alpha(s)\,ds)\,dx = \mu_0 t(u_- + \int_0^t \alpha(s)\,ds).$$

PROOF. Using Theorem 2 (and the curves (11)) we have that $\lim_{\mu \to \mu_0} (u_*^\mu + \int_0^t \alpha(s)ds) = \lim_{\mu \to \mu_0} \frac{dx_1(t)}{dt} = u_- + \int_0^t \alpha(s)ds$ and

$$\lim_{\mu \to \mu_0} \frac{dx_2(t)}{dt} = \lim_{\mu \to \mu_0} \left( u_+ + \frac{\mu}{\rho_+} + \int_0^t \alpha(s)ds \right) = u_- + \int_0^t \alpha(s)ds.$$

On the other hand, we have

$$\lim_{\mu \to \mu_0} \int_{x_1(t)}^{x_2(t)} \rho_*^\mu dx = \mu_0 t,$$

and

$$\lim_{\mu \to \mu_0} \int_{x_1(t)}^{x_2(t)} \rho_*^\mu (u_*^\mu + \int_0^t \alpha(s)ds) dx = (u_- + \int_0^t \alpha(s)ds) \lim_{\mu \to \mu_0} \int_{x_1(t)}^{x_2(t)} \rho_*^\mu dx = \mu_0 t(u_- + \int_0^t \alpha(s)ds).$$

The equality (43) and Lemmas 12 and 13 show that the limit of $\rho$ as $\mu \to \mu_0$ has the same singularity as a weighted Dirac delta function at $x_\delta(t)$ where $\frac{dx_\delta(t)}{dt} = \sigma(t)$. Now, we show the theorem characterizing the limit as $\mu \to \mu_0$.

**Theorem 14.** *Let $u_- < u_+$, $\rho_+(u_- - u_+) = \mu_0 < \mu$, and assume that $(\rho, u)$ is the Riemann solution $J_1 + J_2$ to the problem* (1)-(3) *constructed in Section 2. Then,*

$$\lim_{\mu \to \mu_0} (\rho^\mu(x,t), u^\mu(x,t)) = \begin{cases} (\rho_-, u_- + \int_0^t \alpha(s)ds), & \text{if } x < x_\delta(t), \\ (w(t)\delta(x_\delta(t) - x), \sigma(t)), & \text{if } x = x_\delta(t), \\ (\rho_+, u_+ + \int_0^t \alpha(s)ds), & \text{if } x > x_\delta(t), \end{cases}$$

*in the sense of distributions and where $w(t) = \mu_0 t$, $\sigma(t)$ is given by* (44)*, and $x_\delta(t) = u_- t + \int_0^t \int_0^r \alpha(s)ds\,dr$.*

PROOF. Using the similar variable (35) together the boundary problem (38)-(39), from the first equation, we have a weak formulation

$$\int_{-\infty}^{\infty} \rho(\widetilde{u} - \xi)\phi' d\xi + \int_{-\infty}^{\infty} \rho\phi d\xi = 0 \tag{46}$$

for any $\phi \in C_0^\infty(\mathbb{R})$. The first integral on the left-hand side of (46) can be decomposed into

$$\int_{-\infty}^{\infty} \rho(\widetilde{u} - \xi)\phi' d\xi = \int_{-\infty}^{u_-} \rho(\widetilde{u} - \xi)\phi' d\xi + \phi'd\xi + \int_{u_-}^{u_+ + \mu/\rho_+} \rho(\widetilde{u} - \xi)\phi' d\xi + \int_{u_+ + \mu/\rho_+}^{\infty} \rho(\widetilde{u} - \xi)\phi' d\xi \tag{47}$$

The limit of the sum of the first and last term of (47) is

$$\begin{aligned} \lim_{\mu \to \mu_0} & \left( \int_{-\infty}^{u_-} \rho(\widetilde{u} - \xi)\phi' d\xi + \int_{u_+ + \mu/\rho_+}^{\infty} \rho(\widetilde{u} - \xi)\phi' d\xi \right) \\ &= \lim_{\mu \to \mu_0} \int_{-\infty}^{u_-} \rho_-(u_- - \xi)\phi' d\xi + \lim_{\mu \to \mu_0} \int_{u_+ + \mu/\rho_+}^{\infty} \rho_+(u_+ - \xi)\phi' d\xi \\ &= \rho_- u_- \phi(u_-) - \rho_- \int_{-\infty}^{u_-} \xi\phi' d\xi \rho_+ u_+ \phi(u_-) - \rho_+ \lim_{\mu \to \mu_0} \int_{u_+ + \mu/\rho_+}^{\infty} \xi\phi' d\xi \\ &= \rho_+(u_- - u_+)\phi(u_-) + \int_{-\infty}^{\infty} h(\xi - u_-)\phi d\xi = \mu_0\phi(u_-) + \int_{-\infty}^{\infty} h(\xi - u_-)\phi d\xi \end{aligned} \tag{48}$$

with

$$h(x) = \begin{cases} \rho_-, & \text{if } x < 0, \\ \rho_+, & \text{if } x > 0. \end{cases} \tag{49}$$

The limit of the second term of (47) is

$$
\begin{aligned}
\lim_{\mu \to \mu_0} \int_{u_-}^{u_+ + \mu/\rho_+} \rho(\widetilde{u} - \xi)\phi' d\xi &= \lim_{\mu \to \mu_0} \int_{u_-}^{u_+ + \mu/\rho_+} \rho_*(\widetilde{u}_* - \xi)\phi' d\xi = \lim_{\mu \to \mu_0} \int_{u_-}^{u_+ + \mu/\rho_+} \frac{\mu\rho_+}{\mu + \rho_+(u_+ - u_-)}(u_- - \xi)\phi' d\xi \\
&= -\lim_{\mu \to \mu_0} \frac{\mu\rho_+ u_-}{\mu + \rho_+(u_+ - u_-)} \left( (u_+ + \frac{\mu}{\rho_+})\phi(u_+ + \frac{\mu}{\rho_+}) - u_-\phi(u_-) - \int_{u_-}^{u_+ + \mu/\rho_+} \phi d\xi \right) = 0.
\end{aligned}
$$

(50)

From (48), (50), and (46) we have

$$
\lim_{\mu \to \mu_0} \int_{-\infty}^{\infty} \rho\phi d\xi = \mu_0\phi(u_-) + \int_{-\infty}^{\infty} h(\xi - u_-)\phi d\xi
$$

with $h$ defined by (49). From the previous analysis, for any $\varphi \in C_0^\infty(\mathbb{R} \times (0, \infty))$,

$$
\begin{aligned}
\lim_{\mu \to \mu_0} \int_0^\infty \iint_{\mathbb{R}} \rho(\frac{x - \int_0^t \int_0^r \alpha(s)ds}{t})\varphi(x, t)dxdt \\
&= \lim_{\mu \to \mu_0} \int_0^\infty \iint_{\mathbb{R}} \rho(\xi)\varphi(\xi t + \int_0^t \int_0^r \alpha(s)dsdr, t)d(\xi t + \int_0^t \int_0^r \alpha(s)dsdr)dt \\
&= \lim_{\mu \to \mu_0} \int_0^\infty t \left( \int_{\mathbb{R}} \rho(\xi)\varphi(\xi t + \int_0^t \int_0^r \alpha(s)dsdr, t)d\xi \right) dt \\
&= \int_0^\infty t \left( \mu_0\varphi(u_-t + \int_0^t \int_0^r \alpha(s)dsdr, t) + \int_{-\infty}^\infty h(\xi - u_-)\varphi(\xi t + \int_0^t \int_0^r \alpha(s)dsdr, t)d\xi \right) dt \\
&= \int_0^\infty \mu_0 t\varphi(u_-t + \int_0^t \int_0^r \alpha(s)dsdr, t)dt + \int_0^t t \left( \int_{-\infty}^\infty h(\xi - u_-)\varphi(\xi t + \int_0^t \int_0^r \alpha(s)dsdr, t)d\xi \right) dt \\
&= \int_0^\infty \mu_0 t\varphi(u_-t + \int_0^t \int_0^r \alpha(s)dsdr, t)dt + \int_0^t \int_{-\infty}^\infty h(x - u_-t - \int_0^t \int_0^r \alpha(s)dsdr)\varphi(x, t)dxdt
\end{aligned}
$$

and from Definition 1 we get

$$
\int_0^\infty \mu_0 t\varphi(u_-t + \int_0^t \int_0^r \alpha(s)dsdr, t)dt = <w(\cdot)\delta_L, \varphi(\cdot, \cdot)>
$$

with $w(t) = \mu_0 t$. In addition, using (44) we can obtain the limit of $u$ as $\mu \to \mu_0$ and we omit. Finally, from (44) again, to obtain the curve $x_\delta(t)$ we define $\frac{dx_\delta(t)}{dt} = \lim_{\mu \to \mu_0} \frac{dx_1(t)}{dt} = \sigma(t) = u_- + \int_0^t \alpha(s)ds$ and we have that $x_\delta(t) = u_-t + \int_0^t \int_0^r \alpha(s)dsdr$. It is clear that the entropy condition (13) is valid.

**Theorem 15.** *If $u_- > u_+$, $\mu < \mu_0 = \rho_+(u_- - u_+)$, and assume that $(\rho, u)$ is the delta-shock to the problem (4)-(3) constructed in Section 3. Then, when $\mu \to 0$, $(\rho^\mu, u^\mu)$ converges to a delta-shock wave solution for the pressureless system with a time-dependent Coulomb-like friction term.*

PROOF. As $\mu < \mu_0$, then $(\rho_+, u_+) \in V$ and the Riemann solution to the problem (4)-(3) is the delta-shock given by (27) with functions $u_\delta(t)$, $x_\delta(t)$, and $w(t)$ according to (28) or (29). Observe that as $\mu \to 0$, we have

$$
\lim_{\mu \to 0}(\rho^\mu(x, t), u^\mu(x, t)) = \begin{cases} (\rho_-, u_- + \int_0^t \alpha(s)ds), & \text{if } x < \gamma(t), \\ (\overline{w}(t)\delta(\gamma(t) - x), u_\delta^\rho(t)), & \text{if } x = \gamma(t), \\ (\rho_+, u_+ + \int_0^t \alpha(s)ds), & \text{if } x > \gamma(t), \end{cases}
$$

where $u_\delta^p(t)$, $\gamma(t)$, and $\overline{w}(t)$ are

$$\begin{cases} \lim\limits_{\mu \to 0} u_\delta(t) = \frac{1}{2}\left(u_- + u_+\right) + \int_0^t \alpha(s)ds =: u_\delta^p(t), \\ \lim\limits_{\mu \to 0} x_\delta(t) = \frac{1}{2}\left(u_- + u_+\right)t + \int_0^t \int_0^r \alpha(s)ds dr =: \gamma(t), \qquad \text{if } \rho_- - \rho_+ = 0, \\ \lim\limits_{\mu \to 0} w(t) = \rho_+(u_- - u_+)t =: \overline{w}(t), \end{cases}$$

or

$$\begin{cases} \lim\limits_{\mu \to 0} u_\delta(t) = \frac{\sqrt{\rho_-}u_- + \sqrt{\rho_+}u_+}{\sqrt{\rho_-} + \sqrt{\rho_+}} + \int_0^t \alpha(s)ds =: u_\delta^p(t), \\ \lim\limits_{\mu \to 0} x_\delta(t) = \frac{\sqrt{\rho_-}u_- + \sqrt{\rho_+}u_+}{\sqrt{\rho_-} + \sqrt{\rho_+}}t + \int_0^t \int_0^r \alpha(s)ds dr =: \gamma(t), \qquad \text{if } \rho_- - \rho_+ \neq 0. \\ \lim\limits_{\mu \to 0} w(t) = \sqrt{\rho_- \rho_+}(u_- - u_+)t =: \overline{w}(t), \end{cases}$$

## 5. Comments, Examples, and Extensions

The main goal of this section is to present some comments on the paper and consider two straightforward examples of the Riemann solutions for both systems (1) and (2). Moreover, we discuss about extensions on the subject elaborated upon in this paper.

### 5.1. Comments

In this work, we study the Riemann problem for two systems: the nonsymmetric system of the Keyfitz-Kranzer type system and a pressureless system, both with a time-dependent Coulomb-like friction term. Both systems are important for describing physical models as well as interesting and new mathematical results. Riemann solutions are characterized by contact discontinuities and delta-shock waves for both systems (1) and (2). Moreover, the existence and uniqueness of Riemann solutions, including delta-shock solutions, were obtained for both systems.

In addition, we provide an affirmative response to the question regarding the convergence of Riemann solutions from the nonsymmetric Keyfitz-Kranzer system (1) to those of the pressureless system (2) as $\mu$ tends to 0, which is supported by rigorous mathematical analysis.

### 5.2. Examples

To illustrate the application of our results in the above sections, we present two examples. In these two models, we focus our attention on Riemann solutions involving contact discontinuities and delta-shocks.

#### 5.2.1. Constant coefficient Coulomb-like friction term

In this section, we consider $\alpha(\cdot) \equiv \alpha = const.$. The Riemann solution to the problem (1)-(3) for $u_- < u_+ + \mu/\rho_+$ (using Theorem 2) is

$$(\rho(x,t), u(x,t)) = \begin{cases} (\rho_-, u_- + \alpha t), & \text{if } x < u_- t + \alpha \frac{t^2}{2}, \\ \left(\frac{\mu \rho_+}{\mu + \rho_+(u_+ - u_-)}, u_- + \alpha t\right), & \text{if } u_- t + \alpha \frac{t^2}{2} < x < (u_+ + \frac{\mu}{\rho_+})t + \alpha \frac{t^2}{2}, \\ (\rho_+, u_+ + \alpha t), & \text{if } x > (u_+ + \frac{\mu}{\rho_+})t + \alpha \frac{t^2}{2}, \end{cases}$$

and for $u_- > u_+ + \mu/\rho_+$, by Theorem 6, the solution is

$$(\rho(x,t), u(x,t)) = \begin{cases} (\rho_-, u_- + \alpha t), & \text{if } x < x_\delta(t), \\ (w(t)\delta(x_\delta(t) - x), u_\delta(t)), & \text{if } x = x_\delta(t), \\ (\rho_+, u_+ + \alpha t), & \text{if } x > x_\delta(t), \end{cases}$$

where $u_\delta(t)$, $x_\delta(t)$, and $w(t)$ are

$$\begin{cases} u_\delta(t) = \frac{1}{2}\left(u_- + u_+ + \frac{\mu}{\rho_+}\right) + \alpha t, \\ x_\delta(t) = \frac{1}{2}\left(u_- + u_+ + \frac{\mu}{\rho_+}\right)t + \alpha\frac{t^2}{2}, & \text{if } \rho_- - \rho_+ = 0, \\ w(t) = \rho_+(u_- - u_+)t, \end{cases}$$

or

$$\begin{cases} u_\delta(t) = \frac{[\rho u] - \sqrt{[\rho u]^2 - [\rho][\rho u^2 + \mu u]}}{[\rho]} + \alpha t, \\ x_\delta(t) = \frac{[\rho u] - \sqrt{[\rho u]^2 - [\rho][\rho u^2 + \mu u]}}{[\rho]}t + \alpha\frac{t^2}{2}, & \text{if } \rho_- - \rho_+ \neq 0. \\ w(t) = \sqrt{[\rho u]^2 - [\rho][\rho u^2 + \mu u]}t, \end{cases}$$

On the other hand, the pressureless system (2) has a Riemann solution given by

$$(\rho(x,t), u(x,t)) = \begin{cases} (\rho_-, u_- + \alpha t), & \text{if } x < u_- t + \frac{1}{2}\alpha t^2, \\ (0, \frac{x - \frac{1}{2}\alpha t^2}{t} + \alpha t), & \text{if } u_- t + \frac{1}{2}\alpha t^2 < x < u_+ t + \frac{1}{2}\alpha t^2, \\ (\rho_+, u_+ + \alpha t), & \text{if } x > u_+ t + \frac{1}{2}\alpha t^2, \end{cases}$$

for $u_- < u_+$ (see Theorem 9) and

$$(\rho(x,t), u(x,t)) = \begin{cases} (\rho_-, u_- + \alpha t), & \text{if } x < x_\delta(t), \\ (w(t), \delta(x_\delta(t) - x)), & \text{if } x = x_\delta(t), \\ (\rho_+, u_+ + \alpha t), & \text{if } x > x_\delta(t), \end{cases}$$

for $u_- > u_+$ and where $w(t) = \sqrt{\rho_-\rho_+}(u_- - u_+)t$, $x_\delta(t) = \frac{\sqrt{\rho_-}u_- + \sqrt{\rho_+}u_+}{\sqrt{\rho_-} + \sqrt{\rho_+}} + \frac{1}{2}\alpha t^2$, and $u_\delta(t) = \frac{\sqrt{\rho_-}u_- + \sqrt{\rho_+}u_+}{\sqrt{\rho_-} + \sqrt{\rho_+}} + \alpha t$.

### 5.2.2. Time-gradually-degenerate friction term

In the literature, the time-dependent coefficient $\alpha(t) = -\frac{\theta}{(1+t)^\beta}$ with physical parameters $\theta > 0$ and $\beta \geq 0$ is a *time-gradually-degenerate friction term* [15]. Therefore, the Coulomb-like friction term decayed as $\frac{\theta}{(1+t)^\beta}$ and represents the time-gradually vanishing friction effect. Observe that for the time-gradually-degenerate friction term, we have

$$\int_0^t \frac{\theta}{(1+s)^\beta}ds = \begin{cases} \theta\ln(t+1), & \text{if } \beta = 1, \\ \theta\frac{t}{t+1}, & \text{if } \beta = 2, \\ \frac{\theta}{1-\beta}((t+1)^{1-\beta} - 1), & \text{if } \beta \geq 0, \beta \neq 1, \text{ and } \beta \neq 2, \end{cases} \tag{51}$$

and

$$\int_0^t \int_0^r \frac{\theta}{(1+s)^\beta}ds\,dr = \begin{cases} \theta((t+1)\ln(t+1) - t), & \text{if } \beta = 1, \\ \theta(t - \ln(t+1)), & \text{if } \beta = 2, \\ \frac{\theta}{1-\beta}(\frac{1}{2-\beta}((t+1)^{2-\beta} - 1) - t), & \text{if } \beta \geq 0, \beta \neq 1, \text{ and } \beta \neq 2. \end{cases} \tag{52}$$

To the nonsymmetric Keyfitz-Kranzer (1), according to Theorem 2, when $u_- < u_+ + \mu/\rho_+$, the Riemann solution is given by

$$(\rho(x,t), u(x,t)) = \begin{cases} (\rho_-, u_- + \int_0^t \frac{\theta}{(1+s)^\beta}ds), & \text{if } x < u_- t + \int_0^t \int_0^r \frac{\theta}{(1+s)^\beta}ds\,dr, \\ (\rho_*, u_- + \int_0^t \frac{\theta}{(1+s)^\beta}ds), & \text{if } u_- t + \int_0^t \int_0^r \frac{\theta}{(1+s)^\beta}ds\,dr < x < (u_- + \frac{\mu}{\rho_+})t + \int_0^t \int_0^r \frac{\theta}{(1+s)^\beta}ds\,dr, \\ (\rho_-, u_- + \int_0^t \frac{\theta}{(1+s)^\beta}ds), & \text{if } x > (u_+ + \frac{\mu}{\rho_+})t + \int_0^t \int_0^r \frac{\theta}{(1+s)^\beta}ds\,dr, \end{cases}$$

where $\rho_* = \frac{\mu \rho_+}{\mu + \rho_* (u_+ - u_-)}$, and the integrals $\int_0^t \frac{\theta}{(1+s)^\beta} ds$ and $\int_0^t \int_0^r \frac{\theta}{(1+s)^\beta} ds dr$ are given by (51) and (52), respectively. On the other hand, from Theorem 6, for $u_- > u_+ + \mu/\rho_+$ the Riemann solution is given by the following delta-shock solution:

$$(\rho(x,t), u(x,t)) = \begin{cases} (\rho_-, u_- + \int_0^t \frac{\theta}{(1+s)^\beta} ds), & \text{if } x < x_\delta(t), \\ (w(t)\delta(x_\delta(t) - x), u_\delta(t)), & \text{if } x = x_\delta(t), \\ (\rho_+, u_+ + \int_0^t \frac{\theta}{(1+s)^\beta} ds), & \text{if } x > x_\delta(t), \end{cases}$$

where $u_\delta(t)$, $x_\delta(t)$, and $w(t)$ are

$$\begin{cases} u_\delta(t) = \frac{1}{2}\left(u_- + u_+ + \frac{\mu}{\rho_+}\right) + \int_0^t \frac{\theta}{(1+s)^\beta} ds, \\ x_\delta(t) = \frac{1}{2}\left(u_- + u_+ + \frac{\mu}{\rho_+}\right)t + \int_0^t \int_0^r \frac{\theta}{(1+s)^\beta} ds dr, \qquad \text{if } \rho_- - \rho_+ = 0, \\ w(t) = \rho_+(u_- - u_+)t, \end{cases}$$

or

$$\begin{cases} u_\delta(t) = \frac{[\rho u] - \sqrt{[\rho u]^2 - [\rho][\rho u^2 + \mu u]}}{[\rho]} + \int_0^t \frac{\theta}{(1+s)^\beta} ds, \\ x_\delta(t) = \frac{[\rho u] - \sqrt{[\rho u]^2 - [\rho][\rho u^2 + \mu u]}}{[\rho]} t + \int_0^t \int_0^r \frac{\theta}{(1+s)^\beta} ds dr, \qquad \text{if } \rho_- - \rho_+ \neq 0, \\ w(t) = \sqrt{[\rho u]^2 - [\rho][\rho u^2 + \mu u]} t, \end{cases}$$

and the integrals $\int_0^t \frac{\theta}{(1+s)^\beta} ds$ and $\int_0^t \int_0^r \frac{\theta}{(1+s)^\beta} ds dr$ are given by (51) and (52), respectively.

Now, for the pressureless system (2), from Theorem 9, for $u_- < u_+$ the Riemann solution is given by

$$(\rho(x,t), u(x,t)) = \begin{cases} (\rho_-, u_- + \int_0^t \frac{\theta}{(1+s)^\beta} ds), & \text{if } x < u_- t + \int_0^t \int_0^r \frac{\theta}{(1+s)^\beta} ds dr, \\ (0, \frac{x - \int_0^t \int_0^r \frac{\theta}{(1+s)^\beta} ds dr}{t} + \int_0^t \frac{\theta}{(1+s)^\beta} ds), & \text{if } u_- t + \int_0^t \int_0^r \frac{\theta}{(1+s)^\beta} ds dr < x < u_+ t + \int_0^t \int_0^r \frac{\theta}{(1+s)^\beta} ds dr, \\ (\rho_+, u_+ + \int_0^t \frac{\theta}{(1+s)^\beta} ds), & \text{if } x > u_+ t + \int_0^t \int_0^r \frac{\theta}{(1+s)^\beta} ds dr, \end{cases}$$

where the integrals $\int_0^t \frac{\theta}{(1+s)^\beta} ds$ and $\int_0^t \int_0^r \frac{\theta}{(1+s)^\beta} ds dr$ are given by (51) and (52), respectively. From Theorem 10, for $u_- > u_+$, the Riemann solution is given by

$$(\rho(x,t), u(x,t)) = \begin{cases} (\rho_-, u_- + \int_0^t \frac{\theta}{(1+s)^\beta} ds), & \text{if } x < x_\delta(t), \\ (w(t)\delta(x_\delta(t) - x), u_\delta(t)), & \text{if } x = x_\delta(t), \\ (\rho_+, u_+ + \int_0^t \frac{\theta}{(1+s)^\beta} ds), & \text{if } x > x_\delta(t), \end{cases}$$

where $w(t) = \sqrt{\rho_- \rho_+}(u_- - u_+)t$, $x_\delta(t) = \frac{\sqrt{\rho_-} u_- + \sqrt{\rho_+} u_+}{\sqrt{\rho_-} + \sqrt{\rho_+}} + \int_0^t \int_0^r \frac{\theta}{(1+s)^\beta} ds dr$, and $u_\delta(t) = \frac{\sqrt{\rho_-} u_- + \sqrt{\rho_+} u_+}{\sqrt{\rho_-} + \sqrt{\rho_+}} + \int_0^t \frac{\theta}{(1+s)^\beta} ds$ and the integrals $\int_0^t \frac{\theta}{(1+s)^\beta} ds$ and $\int_0^t \int_0^r \frac{\theta}{(1+s)^\beta} ds dr$ are given by (51) and (52), respectively.

### 5.3. Extensions

One remarks that, the results settled in this paper combined with the ones established by the authors in [9] apply to some correlated versions of the systems (1) and (2), that is, we may consider the following systems:

$$\begin{cases} \rho_t + (\rho u)_x = 0, \\ (\rho u)_t + (\rho u(u + \frac{\mu}{\rho}))_x = \alpha(t)\rho - \sigma(t)\rho u \end{cases}$$

and

$$\begin{cases} \rho_t + (\rho u)_x = 0, \\ (\rho u)_t + (\rho u^2)_x = \alpha(t)\rho - \sigma(t)\rho u, \end{cases} \qquad (53)$$

where $\sigma$ is a nonnegative function, and the term $\sigma(t)\rho u$ represents a time-gradually-degenerate damping. In particular, we may consider the Eulerian droplet model in conservative form, which is (53) with $\alpha(t) = \sigma(t)u_a(t)$, where $\rho$ and $u$ denote the volume fraction and velocity of the particles (droplets), respectively, $u_a$ is a function which depends to time and represents the velocity of the carrier fluid (air), and $\sigma$ is the time-drag coefficient between the carrier fluid and the particles [4, 18].

## CRediT authorship contribution statement

**Richard De la cruz:** Conceptualization, Formal analysis, Writing – original draft, Supervision, Revision.
**Wladimir Neves:** Conceptualization, Formal analysis, Writing – original draft, Supervision, Revision.

## Declaration of competing interest

The authors declare that they have no known competing financial interests or personal relationships that could have appeared to influence the work reported in this paper.

## Data availability

No data was used for the research described in the article.

## Acknowledgements

The author Richard De la cruz acknowledges the support received from Universidad Pedagógica y Tecnológica de Colombia. The author Wladimir Neves has received research grants from CNPq through the grants 313005/2023-0, 406460/2023-0, and also by FAPERJ (Cientista do Nosso Estado) through the grant E-26/201.139/2021.